\documentclass[a4paper]{article}%\documentclass[12pt]{amsart}
\usepackage[utf8]{inputenc}
\usepackage{amsmath}
\usepackage{amssymb}
\usepackage{amsfonts}
\usepackage{amsthm}
\usepackage{mathrsfs} 
\usepackage{bm}
\usepackage{bbm}
\usepackage{tikz}
\usepackage{centernot}
\usepackage{hyperref}
\usepackage[ruled,vlined]{algorithm2e}
\usepackage[margin=1.2in]{geometry}
\newcommand{\Z}{\Bbb{Z}}

\renewcommand{\P}{\Bbb{P}}

\newcommand{\ep}{\epsilon}

\hypersetup{
    colorlinks=true,
    linkcolor=blue,
    filecolor=magenta,
    urlcolor=blue,
    citecolor=blue
}

\makeatletter
\newtheorem*{rep@theorem}{\rep@title}
\newcommand{\newreptheorem}[2]{%
\newenvironment{rep#1}[1]{%
 \def\rep@title{#2 \ref{##1}}%
 \begin{rep@theorem}}%
 {\end{rep@theorem}}}
\makeatother

\newtheorem{thm}{Theorem}
\newreptheorem{thm}{Theorem}

\newtheorem{result}{Result}[section]
\newtheorem{lem}[result]{Lemma}

\newtheorem{cor}[result]{Corollary}

\theoremstyle{definition}
\newtheorem{rmk}[result]{Remark}

\newtheorem*{ack}{Acknowledgements}

\theoremstyle{remark}

\newcommand{\hide}[1]{}

\newcommand{\rough}[1]{}%\textbf{\textcolor{blue}{#1}}}
\definecolor{darkgreen}{RGB}{75,150,75}
\newcommand{\review}[1]{}%\textcolor{darkgreen}{#1}}

\newcommand{\hides}[1]{}%1}

\newcommand{\pub}[1]{}%\textcolor{purple}{#1}}

\title{A short proof that \texorpdfstring{$w(3,k)\ge (1-o(1))k^2$}{w(3,k) > (1-o(1))k\^2}}
\author{Zach Hunter}
\date{\today}

\begin{document}

\maketitle
\begin{abstract}
    Here we present a short proof that the two-color van der Waerden number $w(3,k)$ is bounded from below by $(1-o(1))k^2$. Previous work has already shown that a superpolynomial lower bound holds for $w(3,k)$. However, we believe our result is still is of interest due to our techniques.
\end{abstract}
\section{Introduction}

In this short note, we discuss the off-diagonal van der Waerden numbers $w(3,k)$, defined as follows. For integer $k\ge 3$, we define $w(3,k)$ to be the smallest $N$ such that for any blue-red coloring of $[N]:= \{1,\dots, N\}$, there either is a blue arithmetic progression of length 3, or a red arithmetic progression of legnth $k$. 

Based off of computational data, Graham had conjectured that $w(3,k) = O(k^2)$. Indeed, it was suggested that $w(3,k) \le (1+o(1))k^2$ might hold.

It turns out that Graham's conjecture is false. Indeed, in a breakthrough paper by Green, it was shown that $w(3,k)$ grows superpolynomially \cite{green}. This was subsequently improved by the author \cite{hunter}, giving the current best known lower bound.
\begin{thm}[{\cite[Theorem~1]{hunter}}]\label{best} We have that
\[w(3,k) \ge \exp(c \log^2 k /\log\log k)\]for some absolute constant $c>0$.
\end{thm}\hide{\noindent We note that it is reasonable to conjecture that Theorem~\ref{best} is close to best possible. Indeed, heuristic arguments suggest that $w(3,k)\le \exp(O(\log^2 k))$ (surpassing this requires new insights about constructing sets $S\subset [N]$ without 3-term arithmetic progressions\footnote{Specifically, one would either need to find denser sets $S$ than what is currently known, or get $[N]\setminus S$ to have shorter arithmetic progression than a random set with it's density. Doing either of these would }).}

Here, we present a proof of a much weaker lower bound.
\begin{thm}\label{weak} We have that 
\[w(3,k)\ge (1-o(1))k^2.\]
\end{thm}\noindent We believe Theorem~\ref{weak} is still of interest, for several reasons. First, Theorem~\ref{weak} gives a better lower bound than all arguments prior to the work of \cite{green} (see Remark~\ref{history} for an account of past bounds). Second, we note that our methods give better bounds for small $k$ (see Remark~\ref{atleast} for details). 

Lastly, our proof is very simple and is hopefully instructive. In particular, we use basic group theory to get colorings which have desirable psuedorandom properties. This method appears to be novel, and likely has further applications; for example, forthcoming work of the author will develop these ideas to improve the lower bounds of the diagonal multicolor van der Waerden numbers \cite{hunter2}.

\begin{rmk}\label{history}\textit{Past bounds (prior to the breakthrough of \cite{green}):} Brown, Landman and Robertson proved a lower bound of $w(3,k) \gg k^{2-1/\log\log k}$ with the Lov\'asz Local Lemma \cite{brown}. This was later optimized by Li and Shu to prove $w(3,k) \gg (k/\log k)^2$ \cite{li}. 

Aaronson et al. proved $w(3,k) \gg k^2$ via a slightly involved combinatorial argument \cite{fox}; this result remains unpublished and attained a worse constant than Theorem~\ref{weak}. We note that we were unaware of this result before we wrote an earlier draft of this argument.

Guo and Warnke proved that $w(3,k) \gg n^2/\log n$ by analyzing a random greedy process \cite{guo}, they were unaware of the result of \cite{fox}.

\end{rmk}
\begin{ack}This work was done under the supervision of Ben Green. We are thankful for his comments on the presentation of this note, especially the suggestion to use the language of direct products. We also thank Atticus Stonestrom and Matt Kwan for helpful feedback.

Parts of this note were prepared at IST Austria, we thank them for their hospitality.\end{ack}
\section{Preliminaries}\label{prelim}

\subsection{Notation}\label{not}
Here we collect a few pieces of notation that we will use.

For positive integer $n$, we write $[n]:= \{1,\dots,n\}$.

For positive integer $k$, we refer to (non-trivial) arithmetic progressions of length $k$ as \textit{$k$-AP's}. We say a set of integers $S$ is \textit{$k$-AP-free} if it does not contain any $k$-AP's as subsets.

We extend the above concept to (additive) groups $G$. In particular, give a group $G$, we say $P\subset G$ is a \textit{$k$-AP with respect to $G$} is there exists $g,d\in G, d\neq 0_G$ such that $P = \{g,g+d,\dots,g+(k-1)d\}$. And similarly, we say $S\subset G$ is \textit{$k$-AP-free} if it does not contain any $k$-AP's (with respect to $G$) as subsets.

Lastly, for positive integer $N$, we write $r_3(N)$ to be the maximum cardinality of a $3$-AP-free subset $S\subset [N]$ (which is considered a subset of $\Z$).

\subsection{Some basic lemmas about groups}

For later use, we recall the following well-known lemma about groups.
\begin{lem}\label{direct} Let $n,m$ be coprime integers. We have
\[\Z/nm\Z \cong \Z/n\Z\times \Z/m\Z.\]
\end{lem}\noindent To see that Lemma~\ref{direct} is true, one may confirm that the map $\phi:a+nm\Z \mapsto (a+n\Z,a+m\Z)$ is an isomorphism.

We also note the following key fact. This observation was previously utilized in a paper by Graham \cite{graham} (which studied a quantity related to $w(3,k)$), and also appears in the unpublished manuscript \cite{fox}. 
\begin{lem}\label{friso}Let $S\subset [n]$ be $3$-AP-free (with respect to $\Z$). Then for $m\ge 2n-1$, we have that \[S+m\Z = \{x+m\Z:x\in S\} \subset \Z/m\Z\] is $3$-AP-free with respect to $\Z/m\Z$.
\begin{proof}Set $G = \Z/m\Z$, and let $\varphi:\Z \to G$ be the homomorphism $x\mapsto x + m\Z$. We shall prove $\varphi(S)$ is $3$-AP-free with respect to $G$.

First, we note that a subset $T$ of a group contains a $3$-AP if and only if there are $x,y,z\in T$ with $x\neq y$ and $x+z = 2y$ (indeed, if $P = \{g,g+d,g+2d\}\subset T$ for $d\neq 0$, then one takes $x=g,y=g+d,z=g+2d$; the converse direction is also straight-forward). 

So, we suppose for sake of contradiction that such $x,y,z\in \varphi(S)$ exist. Then, taking their preimages $\overline{x},\overline{y},\overline{z} \in S$, we have that \[\varphi(\overline{x}+\overline{z}) \equiv 2\varphi(\overline{y}) \pmod{m}\]
\[\implies \overline{x}+\overline{z}-2\overline{y}\in \ker(\varphi) = m\Z.\]
\noindent Since $\overline{x},\overline{y},\overline{z} \in S \subset [n]$, we have that $|\overline{x}+\overline{y}-2\overline{z}| \le 2n-2<m$. 

Then, the above implies that $\overline{x}+\overline{z}-2\overline{y} = 0$. Since $\varphi(\overline{x}) = x \neq y = \varphi(\overline{y})$ (and consequently, $\overline{x}\neq \overline{y}$), this should mean that $\{\overline{x},\overline{y},\overline{z}\}\subset S$ form a $3$-AP, contradicting the assumption that $S$ is $3$-AP-free. Thus such $x,y,z\in \varphi(S)$ cannot occur, as desired.\end{proof}
\end{lem}
\begin{rmk}To use jargon from additive combinatorics, Lemma~\ref{friso} was a basic consequence of the fact that the additive sets $[n]$ and $[n]+m\Z$ are ``Freiman $2$-isomorphic''. 
\end{rmk}

\section{Proof of Theorem~1}

We first need the following key lemma.
\begin{lem}\label{shortex} Let $H_1,H_2$ be groups, and consider their direct product
\[G = H_1\times H_2.\]

Suppose we have sets $S = \{x_1,\dots,x_m\} \subset H_1$ and $T_1,\dots,T_m\subset H_2$ that are each $3$-AP-free in their respective groups.

Then,
\[A:= \bigcup_{i=1}^m \{(x_i,y):y\in T_i\}\]is $3$-AP-free with respect to $G$.
\begin{proof}Consider any $g=(g_1,g_2)\in G$ and $d =(d_1,d_2)\in G\setminus \{0_G\}$. 

Suppose for sake of contradiction that $\{g,g+d,g+2d\}\subset A$. Then, we must clearly have
\[\{g_1,g_1+d_1,g_1+2d_1\}\subset S.\]Because $S$ is $3$-AP-free with respect to $H_1$, this means that $d_1 = 0_{H_1}$ and $g_1  \in S$. 

In particular, there must be some $i\in [m]$ such that $g_1 = x_i$, and $\{g,g+d,g+2d\}$ is contained in the coset $\{(x_i,y):y\in H_2\}$. 

By the assumptions that $d \neq 0_G$ and $d_1= 0_{H_1}$, we must have that $d_2\neq 0_{H_2}$. Hence, $\{g_2,g_2+d_2,g_2+2d_2\}$ is a $3$-AP with respect to $H_2$. However, $A\cap \{(x_i,h): h\in H_2\} = \{(x_i,y):y\in T_i\}$, thus $\{g,g+d,g+2d\} \subset A$ would imply $\{g_2,g_2+d_2,g_2+2d_2\}\subset T_i$, contradicting our assumption that $T_i$ is $3$-AP-free (with respect to $H_2$). \hide{It follows that $\{g,g+d,g+2d\}\not\subset A$, }
\end{proof}
\end{lem}

Using Lemma~\ref{shortex}, we will define a random set $A$ which is always $3$-AP-free with desirable properties.
\begin{lem}\label{key}Let $p$ be a prime. Suppose $r_3(\lfloor p/2\rfloor ) = m$.

Then we can randomly construct a set $A\subset [p^2-p]$ that is $3$-AP-free (with respect to $\Z$), so that for each $p$-AP $P\subset [p^2-p]\subset \Z$ we have
\[\P(A\cap P = \emptyset) = (1-m/(p-1))^m.\]Consequently (by a union bound), when $(1-m/(p-1))^m \le p^{-3}$, we have $w(3,p)>p^2-p$.

\begin{rmk}\label{atleast}A well-known infinite $3$-AP-free set is the set of positive numbers $S\subset \Z$ that don't use the digit ``0'' in ternary (so $S = \{1,2,4,5,7,8,13,\dots\}$). This construction was first noted by Erd\H{o}s and Tur\'an in 1936 \cite{erdos}, and was conjectured to give an optimal lower bound for $r_3(N)$ until the work of Salem and Spencer \cite{salem}.

Since $[N] \cap S$ is 3-AP-free and always has cardinality at least $N^{\log 2/\log 3}\ge N^{.63}$, we get a simple lower bound for $r_3(N)$. Using this with Lemma~\ref{key}, we can see that $w(3,p)> p^2-p$ for all primes $p \ge 2^{25}$.
\end{rmk}
\begin{proof}We will first construct a random set $A_0\subset \Z/(p^2-p)\Z$, and obtain $A$ from this set.

To begin, we note that $(p-1)$ and $p$ are coprime, hence $\Z/(p^2-p)\Z \cong \Z/p\Z \times \Z/(p-1)\Z$ by Proposition~\ref{direct}. We write $H_1$ to denote $\Z/p\Z$ and $H_2$ to denote $\Z/(p-1)\Z$.

Since $r_3(\lfloor p/2\rfloor) = m$, there exists some $S  = \{x_1,\dots,x_m\}\subset [\lfloor p/2\rfloor ]$ which is $3$-AP-free. By Lemma~\ref{friso}, we have that $S+p\Z \subset H_1$ and $S+(p-1)\Z\subset H_2$ are both $3$-AP-free sets in their respective groups.

Now, for each $i \in [m]$, we pick a random element $t_i \in H_2$ independently and uniformly at random. We then set $T_i = t_i+S+(p-1)\Z \subset H_2$, which is also $3$-AP-free with respect to $H_2$ (since this property is invariant under translation). Lastly, we take \[A_0:= \bigcup_{i\in [m]} \{(x_i+\Z/p\Z,y):y\in T_i\}.\]By Lemma~\ref{shortex}, we have that $A_0$ is $3$-AP-free (with respect to $\Z/(p^2-p)\Z$).

Next, we consider the surjective homomorphism $\phi:\Z\to \Z/(p^2-p)\Z;n\mapsto n+(p^2-p)\Z$ and write $\pi$ to denote the bijection $\phi|_{[p^2-p]}$. We take $A = \pi^{-1}(A_0)$.

It remains to confirm that $A$ has the desired properties.

Due to $\pi$ being a bijection, and $\phi$ a homomorphism, we see that for any $k$, if $P\subset [p^2-p]$ is a $k$-AP (wrt $\Z$) then $\pi(P)$ is a $k$-AP (wrt $\Z/(p^2-p)\Z$) with $|\pi(P)| = k$. Hence, the assumption that $A_0$ being $3$-AP-free immediately implies that is $A$ is $3$-AP-free as well.

We next observe that for any $p$-AP $P_0 \subset \Z/(p^2-p)\Z$ with $|P_0| = p$, that we must have $P_0 = \{x,x+d,\dots,x+(k-1)d\}$ for some $d = (d_1,d_2)$ where $d_1\neq 0_{H_1}$ (as otherwise, $P_0$ would be contained in a coset of the subgroup $\{(0,y):y\in H_2\}$, which has size $p-1<|P_0|$). 

Since $H_1$ is a cyclic group of prime order, we have that $\{0_{H_1},d_1,\dots,(p-1)d_1\} = H_1$. So for each $i\in [m]$, $P_0$ intersects the coset $\{(x_i+p\Z,y):y\in H_2\}$. Whence, by the independence of our chosen $t_1,\dots,t_m$, we have that $\P(P_0 \cap A_0 = \emptyset) = (1-m/(p-1))^m$. By the discussion above, this will imply the desired analogue holds for $A$, completing our proof.
\end{proof}
\end{lem}
\noindent By the discussion of Remark~\ref{atleast}, we get the following.
\begin{cor}\label{whenprime}We have $w(3,p)>p^2-p$ for all sufficiently large primes $p$.
\end{cor}Finally, recall that for every $\ep >0$, there exists $k_0$ so that when $k>k_0$ there exists a prime $p \in [(1-\ep)k,k]$ (this is a consequence of the prime number theorem). Hence, applying Corollary~\ref{whenprime}, we get Theorem~\ref{weak}.

\hide{\section{Further applications}

An $(r,m)$-set-coloring 

We say a group $G$ is $(k;r,s)$-set-colorable if there exists a set-coloring $c:G\to \binom{[r]}{s}$ without any monochromatic $k$-AP. Given $G$ and integers $m,s$, we define
$\kappa(G;r,s)$ to be the smallest $k$ such that $G$ is $(k;r,s)$-set-colorable.

\begin{lem}\label{setex}Let 
\[0\longrightarrow H \stackrel{f}{\longrightarrow}  G \stackrel{g}{\longrightarrow} G/H\longrightarrow 0\]
be a short exact sequence.

Also, suppose we have set-colorings $C_1:G/H\to \binom{[r]}{r'},C_2:H\to \binom{[r']}{s}$ which both lack monochromatic $k$-APs.

Then there exists a set-coloring $c:G\to \binom{[r]}{s}$ without monochromatic $k$-APs.

\begin{proof} We construct $c$ as follows.

For each $x\in G/H$, fix some $y=y_x \in g^{-1}(x)$, and some bijection $\phi_x:C_1(x)\to [r']$. For $x\in G/H,z\in H$, we set \[c(y_x+f(z)) = \phi_x^{-1}(C_2(z)).\]

Since each $\phi_x$ is a bijection, it is easy to verify that $c$ is a set-coloring taking values in $\binom{[r]}{s}$.

It remains to confirm that $c$ lacks monochromatic $k$-APs. This follows from Lemma~\ref{shortex}, since each color class is 

\end{proof}
\end{lem}
Lemma~\ref{setex} immediately gives the following.
\begin{cor}Let 
\[0\longrightarrow H \stackrel{f}{\longrightarrow}  G \stackrel{g}{\longrightarrow} G/H\longrightarrow 0\]
be a short exact sequence.

For any $r\ge r'\ge s$, we have \[\kappa(G;r,s) \le \max\{\kappa(G/H;r,r'),\kappa(H;r',s)\}.\]

\end{cor}

\begin{lem}

\end{lem}}

\end{document}